\documentclass[leqno,12pt]{article}
\usepackage{latexsym}
 \usepackage{graphicx}
\usepackage{amsmath}
\usepackage{amssymb}
\usepackage{mathtools}
\usepackage{cite}
\usepackage{color}
\usepackage{subfig}
\usepackage[margin=1.3in, top=1.3in, bottom=1.3in]{geometry}

\usepackage[colorlinks,linkcolor=black,citecolor=blue,urlcolor=blue,bookmarks=false,hypertexnames=true]{hyperref} 

\newcommand{\nc}{\newcommand}
\nc{\nt}{\newtheorem}
\nt{thm}{Theorem}[section]
\nt{cor}[thm]{Corollary}
\nt{prop}[thm]{Proposition}
\nt{lem}[thm]{Lemma}
\nt{defn}[thm]{Definition}
\nt{rem}[thm]{Remark}
\nt{exa}[thm]{Example}
\nt{ass}[thm]{Assumption}
\nt{alg}[thm]{Algorithm}
\nt{conj}[thm]{Conjecture}
\nc{\ip}[2]{\mbox{$\langle #1,#2 \rangle$}}
\nc{\pf}{\noindent{\bf Proof\ \ }}
\nc{\finpf}{\hfill{$\Box$}\linespace}
\nc{\linespace}{\vspace{\baselineskip} \noindent}
\nc{\R}{{\mathbf R}}
\nc{\oR}{\overline{\R}}
\nc{\M}{\mathcal M}
\nc{\C}{\mathcal C}

\nc{\Rn}{{\mathbf R}^n}
\nc{\inT}{\mbox{\rm int}\,}
\nc{\cl}{\mbox{\rm cl}\,}

\def\tto{\;{\lower 1pt \hbox{$\rightarrow$}}\kern -12pt
           \hbox{\raise 2.8pt \hbox{$\rightarrow$}}\;}
\newenvironment{myequation}{\setcounter{equation}{\value{thm}}
   \begin{equation}}{\addtocounter{thm}{1}\end{equation}}
\newenvironment{myeqnarray}{\setcounter{equation}{\value{thm}}
    \begin{eqnarray}}{\setcounter{thm}{\value{equation}}\end{eqnarray}}

\nc{\bmye}{\begin{myequation}}
\nc{\emye}{\end{myequation}}

\begin{document}
\title{
The complexity of first-order optimization methods from a metric perspective}
\author{
\and
A.S. Lewis
\thanks{ORIE, Cornell University, Ithaca, NY.
\texttt{people.orie.cornell.edu/aslewis} 
\hspace{2cm} \mbox{~}
Research supported in part by National Science Foundation Grant DMS-2006990.}
\and
Tonghua Tian
\thanks{ORIE, Cornell University, Ithaca, NY.  
\texttt{tt543@cornell.edu}}
}
\date{\today}
\maketitle

\begin{abstract}
A central tool for understanding first-order optimization algorithms is the Kurdyka-{\L}ojasiewicz inequality.  Standard approaches to such methods rely crucially on this inequality to leverage sufficient decrease conditions involving gradients or subgradients.  However, the KL property fundamentally concerns not subgradients but rather ``slope'', a purely metric notion.  By highlighting this view, and avoiding any use of subgradients, we present a simple and concise complexity analysis for first-order optimization algorithms on metric spaces.  This subgradient-free perspective also frames a short and focused proof of the KL property for nonsmooth semi-algebraic functions.
\end{abstract}
\medskip

\medskip


\section{Introduction}
Over the last two decades, beginning with pioneering works such as \cite{AbsMahAnd2005,desc_semi}, the Kurdyka-{\L}ojasiewicz inequality has become a central tool for the convergence analysis of first-order optimization algorithms.  Viewed as a subgradient property \cite{loja,Lewis-Clarke}, several factors --- the profound reasons underlying the inequality's validity for a large class of nonsmooth functions (semi-algebraic, for example), and the sometimes involved arguments applying the inequality to optimization complexity analysis  --- can obscure the essential simplicity of the technique.  However, as remarked in \cite{tailwag}, the inequality is fundamentally a metric space property.  This view can be powerful and illuminating, both in non-Euclidean applications like \cite{blanchet-bolte,hauer-mazon}, and --- the philosophy we pursue here --- in understanding and proving the inequality, and leveraging it in simple complexity analysis.  The original motivation for the  KL property in analytic \cite{loj}, more general smooth \cite{Kurd}, and nonsmooth \cite{loja} settings was to bound the length of {\em continuous}-time descent trajectories, a topic outside our current scope, but fundamental in general metric spaces \cite{ambrosio}.

In this self-contained exposition we present some basic and broadly applicable complexity consequences of the KL property, with short, elementary proofs.  Some ingredients are familiar from scattered sources, but by reorganizing the proof sequence and avoiding any use of subgradients or other traditional nonsmooth analysis, we condense and unify the development.  To begin, after introducing the KL inequality, we immediately prove a rudimentary complexity result about objective value convergence, valid in any metric space.  This new observation, startling in its breadth and simplicity, illustrates our two central claims.
\begin{itemize}
\item
Contemporary interest in optimization algorithms for non-Euclidean settings such as manifolds \cite{absil,boumal2022intromanifolds,zhang-sra} invites a subgradient-free complexity analysis.
\item
A purely metric development of KL-based complexity analysis promotes not just generality, but also brevity and transparency.
\end{itemize}

Thus motivated, we present a fresh proof of a nonsmooth Kurdyka-{\L}ojasiewicz property for semi-algebraic (and more general ``definable'') functions on Euclidean spaces.  By focusing on the property's metric rather than subdifferential nature, we simplify and clarify earlier developments with no essential loss in applicability.

Returning to our development of complexity results, we study sequences of iterates satisfying a fundamental {\em slope-descent} condition, combining a traditional notion of sufficient decrease with a metric slope property.  We show the satisfiability of this condition in metric spaces, under mild conditions, and derive local and global convergence results for the resulting iterates.   As an illustration of our framework, we conclude by analyzing majorization-minimization methods for optimization on Riemannian manifolds and more general geodesic spaces. 

\section{The Kurdyka-{\L}ojasiewicz property}
Our setting is a metric space $X$.  In practice, $X$ may be a Euclidean or Hilbert space, or a Riemannian manifold, but our presentation is purely metric.  Laying aside the semi-algebraic setting for our proof of the Kurdyka-{\L}ojasiewicz inequality in Section \ref{proof}, we make no assumptions akin to finite-dimensionality until Section \ref{MMalgorithms}, and until the final examples we use no inner products:  our techniques are thus entirely primal in nature.

We seek to minimize a lower semicontinuous  function $f \colon X \to \bar\R = [-\infty,+\infty]$ that is proper, meaning it  never takes the value $-\infty$ and has nonempty domain
\[
\mbox{dom}\, f ~=~ \{x \in X : f(x) < +\infty\}.  
\]
In practice, $f$ is ofter regular in some sense --- perhaps smooth, or convex, or some composition thereof --- but this development requires no such assumptions.

Denote the metric in the space $X$ by $d(\cdot,\cdot)$ and the distance function to any point $x \in X$ by $d_x(\cdot)$.  Following \cite{giorgi}, the {\em slope} $|\nabla f|(x)$ of $f$ at any point $x$ with finite value $f(x)$ is the infimum of those values $\epsilon > 0$ for which $x$ locally minimizes the perturbed function $f + \epsilon d_x$.  At points $x$ where $f(x) = +\infty$, the slope is defined to be $+\infty$. The notation is motivated by the following standard example.

\begin{exa}[Submanifolds in Euclidean spaces] \label{submanifolds}
{\rm \mbox{} \\
In the Euclidean space $\R^n$, consider any $\C^{(1)}$-smooth submanifold $\M$, endowed with the Euclidean metric, and any differentiable function $f \colon \M \to \R$.  At any point $x \in \M$, the covariant gradient of $f$ is an element of the subspace of $\R^n$ tangent to $\M$ at $x$, and its norm is the slope of $f$.  In particular, when $\M = \R^n$, we have the formula
\[
|\nabla f|(x) ~=~ |\nabla f(x)|.
\]
}
\end{exa}

\begin{defn}
{\rm
A {\em desingularizer} is a continuous function $\phi \colon \R_+ \to \R_+$ with value $\phi(0) = 0$ and continuous, strictly positive derivative on $(0,+\infty)$. The desingularizer $\phi$ has {\em moderate growth} if its derivative satisfies $\phi'(\tau) = o(\frac{1}{\tau})$ as $\tau \downarrow 0$.
}
\end{defn}

\noindent
The following definition follows the idea of \cite{tailwag}.  The strategy is to use a desingularizer 
to induce a type of ``sharpness'' by rescaling the {\em gap}: 
\[
g(\cdot) ~=~ f(\cdot) - \inf f.
\]

\begin{defn} \label{kl}
{\rm
Consider a lower semicontinuous proper function $f \colon X \to \bar\R$ that is bounded below.  We say that $f$ has the {\em Kurdyka-{\L}ojasiewicz property (for minimization)} if there exists a desingularizer $\phi$ such that all points $x \in X$ with small, strictly positive gap $g(x)$ satisfy $|\nabla(\phi \circ g)|(x) \ge 1$.
}
\end{defn}

\noindent
We can rewrite the basic inequality using a simple chain rule \cite[Lemma 4.1]{aze2017}:
\[
|\nabla(\phi \circ g)|(x)  ~=~ \phi'\big(g(x)\big) \cdot |\nabla g|(x)
\]
whenever $0 < g(x) < +\infty$.  The proof is a simple application of the classical Mean Value Theorem. 

More robust than the notion of slope is the {\em limiting slope}, which, at any point $\bar x \in \mbox{dom}\, f$, is the quantity
\[
|\overline{\nabla} f|(\bar x) ~=~ \liminf_{x \to \bar x \atop f(x) \to f(\bar x)} |\nabla f|(x).
\]
It is easy to verify that the definition of the KL property above could equivalently be written with the limiting slope replacing the slope, as could the subsequent chain rule.  When the metric space $X$ is Euclidean, the following remarkable identity (\cite[Proposition 8.5]{ioffe-variational}, for example) always relates the limiting slope and the subdifferential $\partial f$ \cite{va}:
\[
\mbox{dist}\big(0,\partial f(x)\big) ~=~|\overline{\nabla} f|(x),
\]
leading to an equivalent subdifferential version of the KL property.  We emphasize, however, that we make no use of subgradients in this exposition.  

When the metric space $X$ is complete, the KL property is equivalent to an ``error bound'':   for all points $x$ with small gap $g(x)$ and all values $r$ satisfying $0 < r < g(x)$, the distance from $x$ to the level set $g^{-1}[0,r]$ is no larger than $\phi\big(g(x)\big) - \phi(r)$.  We make no use of this  observation (see \cite[Theorem 2.1]{refId0}), but it offers some compelling intuition for the power of the inequality.

Results in \cite{loja,Lewis-Clarke} guarantee that, on a Euclidean space, any lower semicontinuous semi-algebraic (or, more generally, ``globally subanalytic'') function with a bounded set of minimizers has the KL property for minimization, and some corresponding desingularizer has moderate growth.  Indeed, there exists such a desingularizer of the form of a concave function  $\phi(\tau) = \kappa \tau^{1-\theta}$, where the exponent $\theta$ lies in the interval $[0,1)$ and $\kappa > 0$.  Such desingularizers, for example, feature in the original {\L}ojasiewicz proof of the KL property for real-analytic functions \cite{loj}, and its smooth subanalytic extension in \cite{Kurd}. We will often assume desingularizers to be concave in convergence proofs:  for a discussion, see \cite{attouch-projections}.

The rather involved original arguments for the nonsmooth KL property in \cite{loja,Lewis-Clarke}, along with a subsequent proof (due to Drusvyatskiy) in \cite{ioffe-variational} that follows more closely the path laid out in \cite{Kurd}, all depend heavily on the geometry of subgradients.  In line with our philosophy in this work, and inspired by another recent proof \cite[Theorem 2.9]{valette}, in Section \ref{proof} we derive the nonsmooth result from the smooth version in \cite{Kurd} using just a simple slope-based argument.

\section{Objective value convergence} \label{value_convergence}
We next follow a theme with roots in \cite{error_conv}, taken up in \cite{AbsMahAnd2005}, and promoted in \cite{noll-kl,bolte-pauwels-mm,desc_semi,bolte_complex_kur}.  Many classical optimization algorithms --- the first-order methods considered in \cite{bolte_complex_kur} such as gradient descent, the proximal point method, the proximal gradient method, variable metric and quasi-Newton methods, trust region methods, alternating projections, Gauss-Seidel schemes and so forth, majorization-minimization procedures considered in \cite{bolte-pauwels-mm}, including sequential quadratic programming schemes and the moving balls method, and some bundle methods \cite{noll-kl} --- generate sequences in $X$  with a certain descent property.  Specifically, for some constant $\delta > 0$, successive iterates 
$x,x_+$ always satisfy
\bmye \label{basic-descent}
f(x) - f(x_+) ~\ge~ \delta \big(|\overline{\nabla} f|(x_+)\big)^2.
\emye
Algorithms outside the Euclidean or Hilbert setting may also guarantee the same condition:  as remarked in \cite{subgradient-curves}, a basic example is the proximal point method \cite{bacak,jost}.

\begin{thm}[Value convergence] \label{value}
On a metric space $X$, consider any lower semicontinuous proper function $f \colon X \to \bar\R$ that is bounded below and satisfies the KL property for minimization, and any sequence $(x_k)$ in $X$ satisfying the descent property (\ref{basic-descent}) and with initial value $f(x_0)$ near $\inf f$.  Then the values $f(x_k)$ decrease to $\inf f$.  If, furthermore, the corresponding desingularizer has moderate growth, then the gap satisfies
\[
f(x_k) - \inf f ~=~ o\Big(\frac{1}{k}\Big).
\]
\end{thm}

\pf
Consider the gap $\tau_k = g(x_k) \ge 0$ at iteration $k$.  By assumption, we have
\[
\tau_{k-1} - \tau_k ~\ge~ \delta \big(|\overline{\nabla} f|(x_k)\big)^2,
\]
so $\tau_k$ is nonincreasing.  If $\tau_k = 0$ at any iteration $k$ then there is nothing to prove, so we can suppose that $\tau_k$ is small and nonzero for all $k$.

If the result fails, then $\tau_k$ decreases to some limit $\bar\tau > 0$.  Using the chain rule, the KL property ensures
\[
\phi'(\tau_k) \cdot |\overline{\nabla} f|(x_k) ~\ge~ 1
\]
for all $k$, and hence 
\bmye \label{key}
\tau_{k-1} - \tau_k ~\ge~ \frac{\delta}{\big(\phi'(\tau_k)\big)^2}.
\emye
But as $k \to \infty$, the left-hand side converges to zero, while the right-hand side converges to the strictly positive limit $\delta / \big(\phi'(\bar\tau)\big)^2$, which is a contradiction.

Suppose now that the desingularizer $\phi$ has moderate growth.  Given any constant $\mu > 0$, since 
$0 < \tau_k \downarrow 0$, we know $\tau_k < \mu$ and 
$\tau_k \phi'(\tau_k) < \sqrt{\mu\delta}$ for all large $k$.  In particular, there exists an integer $m$ such that whenever $k \ge m$,
\[
\tau_{k-1} - \tau_k ~>~ \frac{1}{\mu} \tau_k^2,
\]
and hence
\[
\frac{1}{\tau_{k-1}} ~<~ 
\frac{1}{\tau_k} - \frac{1}{\mu + \tau_k} ~<~ 
\frac{1}{\tau_k} - \frac{1}{2\mu}.
\]
Induction shows 
\[
\frac{1}{\tau_k} ~\ge~ \frac{1}{\tau_m} + \frac{k-m}{2\mu}
\]
so
\[
k\tau_k ~\le~ \frac{2\mu\tau_m k}{2\mu + \tau_m (k-m) }.
\]
We deduce
\[
\limsup_{k \to \infty} k\tau_k ~\le~ 2\mu.
\]
Since the constant $\mu > 0$ was arbitrary, the result follows.
\finpf

The $o(\frac{1}{k})$ rate is reminiscent of the behavior of the proximal point method for convex optimization \cite[Theorem 3.1]{guler}.  In the case of the standard desingularizer $\phi(\tau) = \kappa \tau^{1-\theta}$ for $\kappa > 0$ and $0\le \theta < 1$, we can refine this result slightly.  The key inequality (\ref{key}) becomes
\[
\tau_{k-1} - \tau_k ~\ge~ \mu \tau_k^{2\theta}
\]
for some constant $\mu > 0$.
An elementary induction then shows four different scenarios (cf.\ \cite{attouch-bolte}) for the rate at which the gap converges to zero:  finite convergence ($\theta = 0$), superlinear 
($0 < \theta < \frac{1}{2}$), linear ($\theta = \frac{1}{2}$), and finally, for 
$\frac{1}{2} < \theta < 1$,
\[
f(x_k) - \inf f  ~=~
O\left(\Big( \frac{1}{k} \Big)^{\frac{1}{2\theta-1}}\right).
\]
The first case is immediate, and the second and third follow from the inequality 
\[
\frac{\tau_{k-1}}{\tau_k} ~\ge~ 1 + \mu\tau_k^{2\theta-1} 
\]
for all large $k$.  In the second case, the right-hand side grows indefinitely with $k$, and in the third it equals $1+\mu$.  Suppose the last case fails.   For any constant $M \ge 1$ suppose
\[
\tau_k ~>~ M\Big( \frac{1}{k} \Big)^{\frac{1}{2\theta-1}}
\]
for some $k \ge 2$ and yet the same inequality fails when we replace $k$ by $k-1$.  We obtain
\[
M\Big( \frac{1}{k} \Big)^{\frac{1}{2\theta-1}} ~+~ 
\mu M^{2\theta} \Big( \frac{1}{k} \Big)^{\frac{2\theta}{2\theta-1}} < 
\tau_k + \mu \tau_k^{2\theta} ~\le~ \tau_{k-1} ~\le~ M\Big( \frac{1}{k-1} \Big)^{\frac{1}{2\theta-1}}
\]
and then convexity implies a contradiction for sufficiently large $M$:
\[
1 + \frac{\mu M^{2\theta - 1}}{k} ~<~ \Big(1 - \frac{1}{k}\Big)^{-\frac{1}{2\theta-1}}
~\le~ \Big(1 - \frac{2}{k}\Big) + \frac{2}{k}\Big(\frac{1}{2} \Big)^{-\frac{1}{2\theta-1}}
~=~ 1 ~+~ \frac{2(2^\frac{1}{2\theta-1} - 1)}{k}.
\]
The induction follows.

\section{A slope-based proof of the KL property} \label{proof}
Viewed through the lens of the slope, the fact that nonsmooth semi-algebraic functions satisfy a Kurdyka-{\L}ojasiewicz property is a simple consequence of the original smooth results due to {\L}ojasiewicz \cite{loj} and Kurdyka \cite{Kurd}.  Appealing to subgradients, as in the original nonsmooth arguments \cite{loja,Lewis-Clarke}, obscures this simplicity.  In the brief proof we present here, like \cite{valette}, we focus entirely on the slope.

Our exposition in fact applies in a setting much broader than that of semi-algebraic geoemetry, using instead the property known as ``definability in an o-minimal structure''. Readers unfamiliar with definable sets and functions can comfortably focus on the semi-algebraic special case. The accessible reference \cite{Coste-min} gives a concise introduction to definability. 

We fix an o-minimal structure $\mathcal{O}$, with respect to which the various functions and sets we consider are definable. We start with a smooth result (cf.\ \cite[Theorem~11]{Lewis-Clarke}), generalizing \cite[Theorem~1]{Kurd}.

\begin{thm}\label{thm:kl-smooth}
In the Euclidean space $\R^n$, consider any bounded $\C^{(1)}$-smooth definable submanifold $\M$, and any differentiable definable function $f \colon \M \to \R$. Then there exists a concave definable desingularizer $\phi$ and a constant $\rho > 0$ such that
\bmye \label{eq:bounded-manifold}
| \nabla_{\mathcal{M}} (\phi \circ f ) (x) | ~\geq~ 1 \quad \text{whenever} ~ 0 < f(x) < \rho,
\emye
where $\nabla_{\M}$ denotes the covariant gradient on $\M$.
\end{thm}

\pf
Consider the definable function $\psi: (0, \infty) \rightarrow \bar\R$ defined as
\begin{equation*}
\psi(t) ~=~ \inf\{| \nabla_{\mathcal{M}} f(x) | : f(x) = t \}.
\end{equation*}
(As usual, the infimum over the empty set is $+\infty$.) We lose no generality in supposing 
$\liminf_{t \searrow 0} \psi(t) = 0$, since otherwise, for all sufficiently small $\rho > 0$, the desingularizer defined by $\phi(t) = \frac{1}{\rho} t$ for $t \ge 0$ has the desired property.  A standard property of definable functions (\cite[Theorem 2.1]{Coste-min}) then implies $\lim_{t \searrow 0} \psi(t) = 0$.

Consider the definable set
\begin{equation*}
S ~=~ \{(x, t) \in \mathcal{M} \times (0, +\infty) : 
f(x) = t,~ |\nabla_{\mathcal{M}} f(x)| < \psi(t) + t^2\}.
\end{equation*}
Since the manifold $\mathcal{M}$ is bounded, an easy argument shows the existence of a point 
$\Bar{x} \in \R^n$ satisfying $(\Bar{x}, 0) \in \mathrm{cl}\, S$. Further standard properties (see \cite[Theorems 2.1 and 3.2]{Coste-min}) now guarantee that the point $(\Bar{x}, 0)$ is accessible by a smooth curve in $S$.  More precisely, there exists a $\C^{(1)}$-smooth curve
$\gamma \colon (0,1) \to \M$ such that the composition $g = f \circ \gamma$ satisfies
$\big(\gamma(s),g(s)\big) \in S$ for $0 < s < 1$, and as $s \searrow 0$, we have 
$(\gamma(s),g(s)\big) \to (\Bar{x}, 0)$ and $(\gamma'(s),g'(s)\big)$ converges to some limit
$(v,r) \in \R^n \times \R_+$.  We can furthermore assume that the derivative 
$g'$ is either constant or strictly monotone.

Whenever $0<s<1$ we have
\[
g'(s) ~=~ \langle \nabla_{\mathcal{M}} f\big(\gamma (s)\big), \gamma'(s) \rangle
~\le~ 
| \nabla_{\mathcal{M}} f\big(\gamma (s)\big) | \cdot | \gamma'(s) |.
\]
Furthermore, we observe
\begin{eqnarray*}
r &=& \lim_{s \searrow 0} \langle \nabla_{\mathcal{M}} f\big(\gamma (s)\big), \gamma'(s) \rangle \\
&\le& \lim_{s \searrow 0} \Big(| \gamma'(s)| \big( \psi \big(g(s)\big) + g(s)^2 \big)\Big)  ~=~ 0,
\end{eqnarray*}
so in fact $r=0$.  Consequently, the derivative $g'$ must be strictly increasing, and so $g$ is also strictly increasing and satisfies $g(s) < s g'(s)$ for all $s \in (0, 1)$. Take any constant $\rho$ in the interval $\big(0,\lim_{s \nearrow 1} g(s)\big)$ and let $\xi \colon (0,\rho] \to (0,1)$ denote the restriction of the inverse function $g^{-1}$, which is also strictly increasing.

Now consider any point $x \in \mathcal{M}$ satisfying $0 < f(x) < \rho$.  Setting $t=f(x)$ and 
$s = \xi(t)$, then providing $\rho$ is sufficiently small, we have
\begin{eqnarray*}
| \nabla_{\mathcal{M}} (\xi \circ f ) (x) | 
&=&
\xi'\big(f(x)\big) \cdot | \nabla_{\mathcal{M}} f (x) |
~\ge~ 
\xi'(t) \cdot \psi(t) 
~=~ 
\frac{\psi\big(g(s)\big)}{g'(s)} \\
& > & \frac{| \nabla_{\mathcal{M}} f \big(\gamma(s)\big) | - g(s)^2}{g'(s)}
~>~
\frac{1}{|\gamma'(s)|} - s g(s) ~>~\frac{1}{1 + |v|}.
\end{eqnarray*}
Note that the derivative $\xi' = \frac{1}{g' \circ \xi} > 0$ is strictly decreasing, and 
$\lim_{t \searrow 0} \xi (t) = 0$. Consequently, the function $\phi: \R_+ \rightarrow \R_+$ defined by 
\[
\phi(t) ~=~
\left\{
\begin{array}{ll}
0 & (t=0) \\
(1 + |v|) \cdot \xi (t) & (0 < t < \rho) \\
(1 + |v|) \cdot \xi (\rho) & (t \ge \rho)
\end{array}
\right.
\]
has the desired property.
\finpf

Using a simple stratification argument, we arrive at the following result of Kurdyka-{\L}ojasiewicz type for nonsmooth definable functions.

\begin{thm}
In the Euclidean space $\R^n$, consider any bounded definable set $U$, and any definable function $f \colon U \to \R$. Then there exists a concave definable desingularizer $\phi$ and a constant $\rho > 0$ such that
\bmye\label{eq:kl-definable}
|\nabla (\phi \circ f )| (x) \geq 1 \quad 
\text{whenever} ~ 0 < f(x) < \rho.
\emye
\end{thm}

\pf
A standard stratification \cite{LVDB} partitions $U$ into finitely-many bounded $\C^{(1)}$-smooth definable submanifolds $\mathcal{M}_1, \dots, \mathcal{M}_k$ of $\R^n$ such that each restriction $f|_{\mathcal{M}_i}$ is a differentiable definable function. Applying Theorem \ref{thm:kl-smooth} to each restriction $f|_{\mathcal{M}_i}$ gives a concave definable desingularizer $\phi_i$ and a constant 
$\rho > 0$ such that
\begin{equation*}
    | \nabla_{\mathcal{M}_i} (\phi_i \circ f|_{\mathcal{M}_i} ) (x) | \geq 1 \quad 
    \text{whenever} ~ x \in \M_i~ \mbox{with}~ 0 < f(x) < \rho.
\end{equation*}
Applying \cite[Theorem 2.1]{Coste-min} to the collection of functions $\phi_i - \phi_j$ with $i \ne j$, by reducing $\rho$ if necessary, we deduce the existence of a particular index $j$ such that $\phi_{j}' \geq \phi_{i}'$ on the interval $(0,\rho)$ for each index $i$. Take $\phi = \phi_{j}$. Now consider any point $x \in U$ with $0 < f(x) < \rho$.  For some $i$ we have $x \in \M_i$, and then
\begin{eqnarray*}
 |\nabla (\phi \circ f )| (x) 
 &\ge&  |\nabla (\phi_i \circ f )| (x) ~\ge~ |\nabla (\phi_i \circ f|_{\mathcal{M}_i} )| (x) \\
 &=& | \nabla_{\mathcal{M}_i} (\phi_i \circ f|_{\mathcal{M}_i} ) (x) | ~\geq~ 1
\end{eqnarray*}
where the first inequality follows from the inequality $\phi' \geq \phi_i'$ and the chain rule, and the second inequality follows from the definition of the slope.
\finpf

\noindent
Note that when $U$ and $f$ are subanalytic, by the Puiseux expansion we know the desingularizer can take the form $\phi(\tau) = \kappa \tau^{1-\theta}$ for some exponent $\theta \in [0, 1)$ and constant $\kappa > 0$ (see \cite{Kurd}).

The proof we present here is more concise than the original development in \cite{Lewis-Clarke}, which relies heavily on the geometry of subgradients.  This brevity, however, comes with a cost. As we remarked earlier, the KL inequality (\ref{eq:kl-definable}) can be written equivalently using the limiting slope, and hence as the following property:
\begin{quote}
Subgradients of the function $\phi \circ f$ at the point $x$ all have norm at least~$1$.
\end{quote}
By contrast, \cite{Lewis-Clarke} proves a stronger result, namely the same property for {\em Clarke} subgradients.  In particular, when the function $f$ is Lipschitz, this ensures that the property above holds for all {\em convex combinations} of subgradients.  By relying on the weaker, limiting-slope-based KL property, we potentially narrow the scope of our convergence analysis:  the basic descent condition (\ref{basic-descent}) at the iterate $x$ amounts to requiring the existence at the next iterate $x_+$ of a subgradient of norm no larger than $\frac{1}{\sqrt{\delta}}\sqrt{f(x)-f(x_+)}$, whereas using the stronger KL property is more flexible, allowing Clarke subgradients instead.  On the other hand, concrete methods such as  majorization-minimization algorithms do indeed enforce condition (\ref{basic-descent}), so this distinction has no impact.  For the convergence results presented here, we are unaware of any practical impact deriving from the stronger KL property in \cite{Lewis-Clarke}:  the simpler slope-based version proved here seems adequate.

\section{Convergence to a minimizer}
As we have seen, convergence in value follows from the basic descent condition (\ref{basic-descent}).
In practice, this condition often follows from two {\em slope-descent conditions}: 
there exists constants $\alpha > 0$ and $\beta > 0$ such that successive iterates $x,x_+$ always satisfy 
\begin{myeqnarray}
f(x) - f(x_+) & \ge & \alpha d^2(x,x_+),  \label{slope-descent1}\\
d(x,x_+) & \ge & \beta|\overline{\nabla} f|(x_+).   \label{slope-descent2}
\end{myeqnarray}
\!\!We first observe that these conditions, which originated with the seminal analysis of \cite{error_conv}, are always satisfiable. We now assume that the underlying metric space $X$ is complete.

\begin{prop}[Existence of slope descent sequences]
On a complete metric space $X$, consider any lower semicontinuous function $f \colon X \to \bar\R$ that is bounded below.  Then for any point $x \in \mbox{\rm dom}\, f$ and constants $\alpha, \beta > 0$ satisfying $\alpha \beta < \frac{1}{2}$, there exists a point $x_+ \in X$ satisfying the slope-descent conditions (\ref{slope-descent1}) and (\ref{slope-descent2}).
\end{prop} 

\pf
In fact we prove a stronger property, with the limiting slope replaced by the slope.
The function $h = f + \alpha d_x^2$ is also lower semicontinuous and bounded below. We claim that, at all points $y \in X$, its slope satisfies
\[
|\nabla h| (y) ~\geq~ |\nabla f| (y) - 2 \alpha d(x, y).
\]
To prove this, we can restrict attention to a point $y \in \mbox{dom}\, h$ where $|\nabla h| (y)$ is finite. Consider any constants $\epsilon' > \epsilon > |\nabla h| (y)$. Since $y$ locally minimizes the function $h + \epsilon d_y$, we deduce that all points $z \in X$ sufficiently close to $y$ satisfy
\begin{eqnarray*}
f(y)
&\leq&
f(z) + \alpha \big( d(x, z)^2 - d(x, y)^2 \big) + \epsilon d(y, z) \\
&\leq& f(z) + \alpha \big( 2 d(x, y) d(y, z) + d(y, z)^2 \big) + \epsilon d(y, z) \\
&\leq& f(z) + \big( 2 \alpha d(x, y) + \epsilon' \big) d(y, z),
\end{eqnarray*}
Consequently $y$ also locally minimizes the function 
\[
f + \big( 2 \alpha d(x, y) + \epsilon' \big) d_y 
\]
and we deduce $|\nabla f|(y) \le 2 \alpha d(x, y) + \epsilon'$.  
Letting $\epsilon' \downarrow |\nabla h| (y)$ proves the claim.

If $x$ minimizes the function $h$, then $|\nabla h|(x) = 0$ and hence $|\nabla f|(x) = 0$, so we can choose $x_+ = x$. Otherwise, define  
$t = \frac{1}{2}(h(x) + \inf h)$ and note that the distance $\delta$ from $x$ to the level set 
$\{y \in X : h(y) \leq t\}$ is strictly positive. The Ekeland principle \cite{var_princ} implies the existence of a point $x_+ \in X$ satisfying both the inequalities $h(x_+) \le t$ and $|\nabla h|(x_+) \le (\frac{1}{\beta} - 2 \alpha) \delta$. We deduce  $d(x,x_+) \ge \delta$ and 
\[
|\nabla f|(x_+) ~\leq~ |\nabla h|(x_+) + 2 \alpha d(x,x_+) ~\leq~ \frac{1}{\beta} d(x,x_+),
\]
so inequality (\ref{slope-descent2}) follows.  Inequality (\ref{slope-descent1}) follows from 
$h(x_+) \le t < h(x)$.
\finpf

The argument above shows in particular that any {\em proximal point} $x_+$ (a point that minimizes the function $f + \alpha d_x^2$) must satisfy the slope descent conditions (\ref{slope-descent1}) and (\ref{slope-descent2}) assuming $\alpha\beta \le  \frac{1}{2}$.  Existence (and uniqueness) of proximal points for lower semicontinuous convex functions $f$ is guaranteed on a Hilbert space $X$ by Minty's Theorem \cite{minty}, and also in more general metric spaces \cite{jost,bacak}.  However, in general, unless the metric space $X$ is {\em proper} (meaning that all closed balls are compact), proximal points may not exist.

The following result and its proof condense arguments in \cite{desc_semi}.  We assume concavity of the desingularizer, and use lower semicontinuity of the objective, but rely on no further continuity properties.

\begin{thm}[Iterate convergence] \label{iterate}
On a complete metric space $X$, consider any lower semicontinuous function $f \colon X \to \bar\R$ that is bounded below and satisfies the KL property for minimization, with a concave desingularizer $\phi$, and any sequence $(x_k)$ in $X$ satisfying the slope-descent conditions {\rm (\ref{slope-descent1})} and {\rm (\ref{slope-descent2})}   with initial value $f(x_0)$ near $\inf f$.  Then, $x_k$ converges to a minimizer of $f$.
\end{thm} 

\pf
We can assume that no iterate minimizes $f$, since otherwise the first inequality would ensure that the sequence is constant thereafter.  The objective values $f(x_k)$ decrease monotonically, so are finite, and the KL inequality ensures that the gap $g$ satisfies
\[
1 ~\le~ |\overline{\nabla}(\phi \circ g)|(x_k)  ~=~ \phi'\big(g(x_k)\big) \cdot \overline{|\nabla g|}(x_k).
\]
We deduce that each limiting slope $|\overline{\nabla} g|(x_k) = |\overline{\nabla} f|(x_k)$ is nonzero, and so the distances 
$d_k = d(x_{k-1},x_k)$ are nonzero for all $k>0$.

Using the concavity of the desingularizer $\phi$, we observe, for all $k>0$,
\begin{eqnarray*}
\phi\big(g(x_k)\big) - \phi\big(g(x_{k+1})\big) 
& \ge & 
\phi'\big(g(x_k)\big) \big(g(x_k) - g(x_{k+1})\big) \\
& \ge & 
\frac{f(x_k) - f(x_{k+1})}{|\overline{\nabla} f|(x_k)} ~\ge~ 
\frac{\alpha d_{k+1}^2}{\beta d_k},
\end{eqnarray*}
so
\[
\frac{\beta}{\alpha}\Big( \phi\big(g(x_k)\big) - \phi\big(g(x_{k+1})\big) \Big)
~ \ge ~
\frac{2d_k d_{k+1} - d_k^2}{d_k}
~ \ge ~
2d_{k+1} - d_k.
\]
In terms of the numbers
\[
\lambda_k ~=~ \frac{\beta}{\alpha}\phi\big(g(x_k)\big) + d_k ~>~ 0 \qquad \mbox{for}~ k>0,
\]
we have proved 
\[
d(x_k,x_{k+1}) ~\le~ \lambda_k - \lambda_{k+1}.
\]
It follows that the sequence $(x_k)$ is Cauchy, and hence converges to some limit $x^*$.  But Theorem \ref{value} shows that $f(x_k)$ decreases to $\inf f$, so lower semicontinuity now ensures that $x^*$ minimizes $f$.
\finpf

In general, we know nothing about how fast the error $d(x_k,x^*)$ converges to zero, even in the canonical case of the proximal point method.  However, with some knowledge of the desingularizer $\phi$, the proof shows more.  Following \cite{attouch-bolte,desc_semi}, we have
\[
d(x_k,x^*) ~\le~ \sum_{j=k}^{m-1} d(x_j,x_{j+1}) ~+~ d(x_m,x^*) ~\le~ \lambda_k - \lambda_m + d(x_m,x^*).
\]
Letting $m \to \infty$ shows
\[
d(x_k,x^*) ~ \le ~ \lambda_k ~=~ \frac{\beta}{\alpha}\phi\big(g(x_k)\big) + d_k
~\le~ \frac{\beta}{\alpha}\phi\big(g(x_k)\big) + \sqrt{\frac{g(x_{k-1})}{\alpha}},
\]
and complexity estimates follow. Following the discussion in the earlier section, in the case of the standard desingularizer $\phi(\tau) = \kappa \tau^{1-\theta}$ for $\kappa > 0$ and $0\le \theta < 1$, 
we obtain four different scenarios in terms of the convergence rate: finite ($\theta = 0$), superlinear ($0<\theta<\frac{1}{2}$), linear ($\theta = \frac{1}{2}$), and, for the case
$\frac{1}{2} < \theta < 1$, the scenario
\[
d(x_k,x^*) ~=~ O \left( \Big( \frac{1}{k} \Big)^{\frac{1-\theta}{2\theta-1}}\right).
\]

\section{Convergence to a critical point}
The previous convergence results rely on the initial value of the objective function $f$ being close to its infimum.  Without that assumption, in the absence of any global convex-like assumptions on $f$,
we instead resort to seeking a point $x^* \in X$ that is {\em critical}, meaning that the limiting slope
$|\overline{\nabla} f|(x^*)$ is zero.  For smooth functions on manifolds, for example, this idea coincides with the classical notion.  Clearly, any local minimizer is a critical point.

In this section, we present conditions guaranteeing a {\em global} convergence property:  any sequence satisfying the slope-descent conditions (\ref{slope-descent1}) and (\ref{slope-descent2}) must converge to a critical point.  The approach essentially follows \cite{desc_semi}, in particular relying first on the following property.

\begin{defn} \label{kl-function}
{\rm
Consider a lower semicontinuous proper function $f \colon X \to \bar\R$.  We call $f$ a {\em KL-function} if, for all points 
$\bar x \in X$ with finite value, there exists a concave desingularizer $\phi$ and a constant $\rho > 0$ such that 
\bmye \label{implication}
d(x,\bar x) < \rho \quad \mbox{and} \quad f(x) - \bar f \in (0,\rho) 
\quad \Rightarrow \quad
|\nabla\big(\phi \circ (f - \bar f)\big)|(x) ~\ge~ 1,
\emye
where $\bar f = f(\bar x)$.
}
\end{defn}

Some remarks about this definition are in order.  First, the desired property is automatic at points $\bar x$ that are not critical, as is easy to verify, so the notion really concerns critical points.   Secondly, we note that all lower semicontinuous semi-algebraic (or globally subanalytic) functions are KL-functions \cite{loja,Lewis-Clarke,valette}.  Finally, as before, we can equivalently replace the slope in the definition by the limiting slope.

From Section \ref{proof}, we deduce the following result.

\begin{cor}
All lower semicontinuous semi-algebraic (or, more generally, subanalytic) functions are KL-functions.
\end{cor}

An important feature of any KL-function $f$ is the possibility of fixing a {\em uniform} desingularizer 
over any compact set $\bar X \subset X$ on which $f$ is constant.  More precisely, there exists a desingularizer $\phi$ and constant $\rho > 0$ such that the KL property (\ref{implication}) holds (in terms of both the slope and limiting slope) for 
all points $\bar x \in \bar X$.  The proof is a routine compactness argument: using the desingularizers 
$\phi_{\bar x}$ and constants $\rho_{\bar x}$ for each point $\bar x$, the corresponding neighborhoods cover $\bar X$, so we can select a subcover corresponding to a finite set of points $\{x_i\}$, and then set the constant $\rho = \min_i \rho_{x_i}$ and the desingularizer $\phi = \sum_i \phi_{x_i}$.  As a consequence of this uniformity property, any KL function with a compact set of minimizers has the Kurdyka-{\L}ojasiewicz property for minimization, in the sense of Definition \ref{kl}.

We would not be exaggerating much to assert that all optimization objectives in practice are KL-functions.  The second condition we assume is less universal, but still common (cf.\ \cite[Theorem 2.3.1]{ambrosio} and \cite{curves}).

\begin{defn} 
{\rm
A proper function $f \colon X \to \bar\R$ is {\em continuous on slope-bounded sets} if, for all points  $\bar x \in \mbox{dom}\,f$ and sequences $x_k \to \bar x$, if the values $f(x_k)$ and slopes $|\nabla f|(x_k)$ are uniformly bounded, then $f(x_k) \to f(\bar x)$.
}
\end{defn}

\noindent
Clearly this property holds if $f$ is continuous, or more generally, continuous on its domain, but it also holds for lower semicontinuous functions on a Euclidean space that are convex or, more generally, ``subdifferentially continuous'' \cite{va}.

\begin{thm}[Global convergence] \label{global}
On a complete metric space $X$, consider any KL-function $f \colon X \to \bar\R$ that is bounded below and continuous on slope-bounded sets.  Then any relatively compact sequence $(x_k)$ in $X$ satisfying the slope-descent conditions {\rm (\ref{slope-descent1})} and {\rm (\ref{slope-descent2})} must converge to a critical point of $f$.
\end{thm} 

\pf
The value $\mu = \inf_k f(x_k)$ is finite.  Note  $f(x_k) - f(x_{k+1}) \to 0$, because the values $f(x_k)$ decrease monotonically to 
$\mu$.  The slope-descent conditions (\ref{slope-descent1}) and (\ref{slope-descent2}) now imply
\bmye \label{null}
d(x_k,x_{k+1}) \to 0 \qquad \mbox{and} \qquad |\overline{\nabla} f|(x_k) \to 0.
\emye
Furthermore, if $f(x_k) = \mu$ ever holds, then $|\overline{\nabla} f|(x_k) = 0$ and $x_j = x_k$ for all $j \ge k$, so there is nothing to prove.  Thus we can assume $f(x_k) > \mu$ for all $k$.

The set $\bar X$ of the limit points of the sequence $(x_k)$ is closed, and is contained in the closure of the sequence, which by assumption is compact.  Hence $\bar X$ is nonempty and compact.  An elementary argument implies 
$\mbox{dist}_{\bar X}(x_k) \to 0$.

Consider any subsequence $K \subset {\mathbf N}$ for which the sequence $(x_k)_{k \in K}$ converges.  Denote the limit by $\bar x$.  As $k \to \infty$ in $K$, we know $f(x_k) \to \mu$ and 
$|\overline{\nabla} f|(x_k) \to 0$, so by the definition of the limiting slope, there exist corresponding points $x'_k$ for $k \in K$ such that $x'_k \to \bar x$, $f(x'_k) \to \mu$, and 
$|\nabla f|(x'_k) \to 0$ as $k \to \infty$ in $K$.  Thus $\bar x$ is critical, and continuity on slope-bounded sets implies $f(\bar x) = \mu$.  Thus $f$ is constant on the set $\bar X$, which furthermore consists entirely of critical points.

Following the discussion of uniformity after Definition \ref{kl-function}, since the set $\bar X$ is compact, there exists a concave desingularizer $\phi$ and a constant $\rho > 0$ such that property (\ref{implication}) holds with 
$\bar f = \mu$ for all points $\bar x \in \bar X$.  For all large $k$, we know 
\[
\mbox{dist}_{\bar X}(x_k) < \rho \qquad \mbox{and} \qquad f(x_k) - \mu \in (0,\rho),
\]
so
\[
|\overline{\nabla}(\phi \circ (f - \mu)|(x_k) ~\ge~ 1.
\]
Following exactly the same argument as in the proof of Theorem \ref{iterate}, we deduce that the whole sequence $(x_k)$ converges, and the limit is critical since it lies in $\bar X$.
\finpf

\section{Majorization-minimization methods} \label{MMalgorithms}
To illustrate the philosophy of this work, we present a purely metric approach to the convergence analysis of majorization-minimization (MM) algorithms, conceptually simply iterative techniques for minimizing an objective function $f \colon X \to \R$ over a set $X$. At each point $x \in X$ we associate a model subset $D(x)$ of $X$ that contains $x$, and a {\em majorizing} model function $h_x$:  that is, $h_x \ge f$ with equality at the point $x$.  The simplest examples have $D(x)$ identically equal to $X$.  We then consider the corresponding 
{\em envelope} $F \colon X \to \bar\R$ and set-valued {\em iteration map} $p \colon X \tto X$ defined at each point $x \in X$ by
\[
F(x) = \inf_{D(x)} h_x \qquad \mbox{and} \qquad p(x) ~=~ \mbox{arg}\!\min_{D(x)} h_x
\]
and {\em MM sequences}, in which successive iterates $x,x_+$ always satisfy $x_+ \in p(x)$. 
Majorization-minimization algorithms have been widely used in practice, originally in the statistics literature \cite{dempster}.

A fundamental initial question is whether an MM sequence $(x_k)$ is a {\em minimizing sequence}  for the objective $f$, by which we simply mean $f(x_k) \to \inf f$.  We first collect some elementary observations.

\begin{prop} \label{observation}
Consider an MM method for an objective $f$ over a set $X$.  Then $f$ majorizes the corresponding envelope $F$, and $\inf f = \inf F$.  An MM sequence is a minimizing sequence for $f$ if and only if it is a minimizing sequence for $F$. Furthermore, the iteration map $p$ satisfies
\[
p(\mbox{\rm argmin}\,F) ~\subset~ \mbox{\rm argmin}\,f ~\subset~ \mbox{\rm argmin}\,F, 
\]
with equality if $p$ is everywhere nonempty-valued and the model satisfies 
$h_x(y) > f(y)$ for all distinct points $x \in X$ and $y \in D(x)$.
\end{prop}

\pf
All $x \in X$ satisfy $F(x) \le h_x(x) = f(x)$, so majorization follows.  We also have
\[
\inf_{y \in X} F(y) \le \inf_{y \in X} f(y) \le \inf_{y \in D(x)} h_x(y) = \inf_x F(x),
\]
so $f$ and $F$ have the same infimum and furthermore $\mbox{\rm argmin}\,f \subset \mbox{\rm argmin}\,F$.
To deduce the claim about MM sequences $(x_k)$, notice $f(x_k) \ge F(x_k) \ge f(x_{k+1})$, because 
$F(x_k) = h_{x_k}(x_{k+1})$.   Finally, consider any minimizer $x$ for $F$.  If there exists a point 
$y \in p(x)$, then we have
\[
\inf F = F(x) = h_x(y) \ge f(y) \ge \inf f = \inf F.
\]
Hence equality holds throughout, so $y$ minimizes $f$. If we have the extra strict domination condition, then the above relationship implies $y=x$.
\finpf

\noindent
The final claim may fail if the iteration map is possibly empty-valued.  An example is the objective $f(x) = e^x$ for $x \in X=\R$, with model subset $D(x) = \R$ for all $x$ and model function 
$h_x(y) = (1+|x-y|)e^y$ for $y \in \R$.  The envelope $F$ is identically zero and the iteration map $p$ is always empty-valued.

The proposition above motivates our subsequent analysis, an approach inspired by \cite{bolte-pauwels-mm} but once again avoiding any variational-analytic tools such as subgradients or normal cones.  Rather than applying convergence results like those in Section \ref{value_convergence} directly to the objective $f$, Proposition \ref{observation} opens up the possibility of working instead with the envelope $F$.

We return to the setting where $X$ is a metric space, and consider bivariate functions 
$h \colon X^2 \to \R$ that are locally Lipschitz, or equivalently, locally Lipschitz with respect to each variable separately.
The novelty in our analysis revolves around the following strengthening of the Lipschitz property with respect to the second variable, an apparently new notion.

\begin{defn} \label{approximator}
{\rm
Given a metric space $X$, an {\em approximator} is a locally Lipschitz function $h \colon X^2 \to \R$ such that for each point $\bar x \in X$ there exists a constant $\gamma \ge 0$ so all points $x,y,z \in X$ near $\bar x$ satisfy the inequality
\bmye \label{approx}
|h(x,y) - h(z,y)| ~\le~ \gamma d(x,z)\big(d(x,y) + d(z,y)\big).
\emye
}
\end{defn}

Here are some simple examples.

\begin{exa}[Squared distance] \label{squared}
{\rm
The function defined by $h(x,y) = d^2(x,y)$ for $x,y \in X$ is an approximator, being locally Lipschitz and satisfying
\begin{eqnarray*}
|h(x,y) - h(z,y)| & = & |d^2(x,y) - d^2(z,y)| ~=~ |d(x,y) - d(z,y)|(d(x,y) + d(z,y)| \\
& \le & d(x,z) \big((d(x,y) + d(z,y)\big)
\end{eqnarray*}
for all $z \in X$.  
}
\end{exa}

\begin{exa}[Threshhold] \label{threshhold}
{\rm
For any radius $\rho > 0$, the function defined by
\[
h(x,y) ~=~ \max\{d(x,y)-\rho,0\} ~=~ \mbox{dist}_{B_\rho(x)}(y) \qquad (x,y \in X)
\]
is an approximator.  Being a maximum of Lipschitz functions, it too is Lipschitz.  Furthermore, for any point $z \in X$, if 
$d(x,y) \ge \rho$ and $d(z,y) \ge \rho$, then
\[
h(x,y) - h(z,y) ~=~ d(x,y) - d(z,y) ~\le~ d(x,z) ~\le~ \frac{1}{2\rho} d(x,z)\big(d(x,y) + d(z,y)\big),
\]
while if $d(x,y) \ge \rho \ge d(z,y)$, then
\begin{eqnarray*}
h(x,y) - h(z,y) & = & d(x,y) - \rho ~\le~ d(x,y) - d(z,y) ~\le~ d(x,z) \\
 & \le & \frac{1}{\rho} d(x,z)\big(d(x,y) + d(z,y)\big).
\end{eqnarray*}
We deduce
\[
h(x,y) - h(z,y)  ~\le~ \frac{1}{\rho} d(x,z)\big(d(x,y) + d(z,y)\big),
\]
the case $d(x,y) \le \rho$ being trivial.  The desired property now follows by symmetry.
}
\end{exa}

\noindent
The following example is particularly useful.

\begin{exa}[Linear approximation] \label{linear}
{\rm
In the Euclidean space $X = \R^n$, consider a function $g \colon \R^n \to \R$ whose gradient $\nabla g$ is $L$-Lipschitz.  Then the function defined from the linear approximation by
\[
h(x,y) ~=~ g(x) + \ip{\nabla g(x)}{y-x} \qquad (x,y \in X)
\]
is an approximator.
}
\end{exa}

\noindent
The following observation is simple to verify.

\begin{prop} \label{vector}
Sums and scalar multiples of approximators are also approximators.
\end{prop}

Our development relies on some strong structural assumptions about the underlying feasible region $X$.  We assume in particular that, as a metric space, $X$ is {\em proper}, meaning that all closed balls are compact.  In particular, $X$ is therefore complete and locally compact.  

\begin{prop}[Uniform approximators]
Consider a proper metric space $X$ and a function $h \colon X^2 \to \R$.  Then $h$ is an approximator if and only if for all bounded subsets $X' \subset X$, the function $h$ is Lipschitz on $(X')^2$ and  there exists a constant $\gamma > 0$ such that the approximator inequality (\ref{approx}) holds uniformly for all points $x,y,z \in X'$.
\end{prop}

\pf
The given property clearly implies that $h$ is an approximator, so we focus on the converse.  Assuming that $h$ is an approximator, the uniform Lipschitz property on bounded sets is a standard compactness argument, and the following proof of inequality (\ref{approx}) is similar.  If the property fails, then there exist sequences 
$(x_k)$, $(y_k)$ and $(z_k)$ in the bounded set $X'$ and constants $\gamma_k \to +\infty$ satisfying, for all $k=1,2,3,\ldots$, the inequality
\bmye \label{contradiction}
|h(x_k,y_k) - h(z_k,y_k)| ~>~ \gamma_k d(x_k,z_k)\big(d(x_k,y_k) + d(z_k,y_k)\big).
\emye
Since $X$ is proper, we can take subsequences and suppose $x_k \to \bar x$, $y_k \to \bar y$, and $z_k \to \bar z$.  Denoting by $\lambda$ a Lipschitz constant for $h$ on $X' \times X'$, we deduce
\[
\lambda ~>~ \gamma_k \big(d(x_k,y_k) + d(z_k,y_k)\big),
\]
so both $d(x_k,y_k)$ and $d(z_k,y_k)$ converge to zero.  Consequently we have $\bar x = \bar y = \bar z$, and so inequality (\ref{contradiction}) contradicts the definition of an approximator.
\finpf

\begin{lem}[Distance functions as approximators] \label{lem:approximator}
Consider a proper metric space $X$ and a set-valued mapping $D \colon X \tto X$ such that for each point $x \in X$, the image $D(x)$ is closed and contains $x$.  Then the function defined by 
\[
h(x,y) ~=~ \mbox{\rm dist}_{D(x)}(y)
\]
is an approximator if and only if $h$ is locally Lipschitz and for each point 
$\bar{x} \in X$ there exists a constant $\gamma \geq 0$ such that all points $x,y,z \in X$ near $\bar x$ with $y \in D(z)$ satisfy the inequality
\[
\mbox{\rm dist}_{D(x)}(y) ~\le~ \gamma d(x,z)\big(d(x,y) + d(z,y)\big).
\]
\end{lem}

\pf
The direct implication is immediate, so we focus on the converse. For any point $x \in X$, the function $H(x,\cdot)$, being a distance function, is $1$-Lipschitz. 
Turning to the second variable and the approximator inequality, first note that for all points $x, y, z \in X$, choosing any $u \in \mbox{\rm Proj}_{D(z)}(y)$ and $w \in \mbox{\rm Proj}_{D(x)}(u)$, we have
\[
\mbox{\rm dist}_{D(x)}(y) - \mbox{\rm dist}_{D(z)}(y)
~=~ 
\mbox{\rm dist}_{D(x)}(y) - d(u,y) ~\le~ d(w,y) - d(u,y) ~\le~ d(w, u) 
\]
so
\bmye \label{dist-approximator}
\mbox{\rm dist}_{D(x)}(y) - \mbox{\rm dist}_{D(z)}(y) ~\le~ \mbox{\rm dist}_{D(x)}(u). 
\emye
Since $z \in D(z)$, we deduce $d(u,y) \le d(z,y)$. So when $x, y, z$ are all near $\bar{x}$, $u$ is also near $\bar x$, and our assumption implies
\begin{eqnarray*}
\mbox{\rm dist}_{D(x)}(y) - \mbox{\rm dist}_{D(z)}(y)
&\le& \mbox{\rm dist}_{D(x)}(u) 
~ \le ~ \gamma d(x,z)\big(d(x,u) + d(z,u)\big) \\
& \le & \gamma d(x,z)\big(d(x,y) + d(u,y) + d(z,y) + d(u,y)\big) \\
& \le & \gamma d(x,z)\big(d(x,y) + 3d(z,y)\big) \\
& \le & 3\gamma d(x,z)\big(d(x,y) + d(z,y)\big).
\end{eqnarray*}
The approximator property (\ref{approx}) follows by symmetry, with the constant $3\gamma$.
\finpf

\noindent
Example \ref{threshhold} illustrates exactly this kind of approximator.

\begin{rem}
{\rm
The Lipschitz continuity of the function $h$ defined in Lemma \ref{lem:approximator} amounts to {\em pseudo-Lipschitz continuity} for the mapping $D$:  see \cite{lip_multi}. 
}
\end{rem}

In the development that follows, in addition to assuming that the feasible region $X$ is a proper metric space, we suppose that it is a closed subset of a possibly larger metric space $\bar X$, with the inherited metric.  The role of $\bar X$ is purely to deploy some convexity arguments.  We use the idea of a {\em geodesic} (of unit speed), which is an isometry from a compact interval of $\R$ into $\bar X$;  the  {\em endpoints} and {\em midpoint} are the images of the endpoints of the interval and its midpoint respectively.  In our discussions, we often identify a geodesic with its image. We call $\bar X$ a {\em geodesic space} if every pair of points are the endpoints of a geodesic.  A subset $D$ is {\em geodesically convex} if $D$ contains all geodesics whose endpoints lie in $D$.  In that case, a function $g \colon D \to \R$ is {\em geodesically convex} if its composition with every geodesic in $D$ is convex, and is {\em geodesically $\mu$-strongly convex} (for $\mu > 0$) if its composition with every geodesic in $D$ is $\mu$-strongly convex.  

A geodesic space $\bar X$ is a CAT(0) {\em space} when, for all points $x \in \bar X$, the function $\frac{1}{2}d^2(x,\cdot)$ is geodesically 1-strongly convex.  This latter property is crucial for the proximal point iteration \cite{jost,bacak}.  In particular, it implies that balls in $\bar X$ are geodesically convex.  Complete, locally compact geodesic spaces are proper as a consequence of a suitable version of the Hopf-Rinow Theorem \cite[Theorem 7.2]{myers} and \cite[Proposition I.3.7]{bridson}.  Examples of such spaces that are also CAT(0) include all finite-dimensional Riemannian manifolds that are complete, simply connected, and having nonpositive sectional curvature. 

For convenience, we collect our assumptions below.

\begin{ass} \label{overall}
The feasible region $X$ is a proper metric subspace of a geodesic space, the objective function $f \colon X \to \R$ is bounded below, and the optimization problem $\inf_X f$ admits the following majorization-minimization tools.
 \begin{itemize}
\item
An {\bf approximating mapping} $D \colon X \tto X$ that assigns to each feasible point \mbox{$x \in X$} a closed and geodesically convex subset $D(x)$ that contains $x$.
\item
A constant $\mu \ge 0$ and a {\bf modeling function} $h \colon X^2 \to \R$ assigns to each $x \in X$
a geodesically $\mu$-strongly convex function $h_x \colon D(x) \to \R$ satisfying
\[
h(x,y) = h_x(y) \ge f(y) \quad \mbox{for}~ y \in D(x), \quad \mbox{with equality if}~ y=x.
\]
\end{itemize}
Furthermore, both $h$ and the function defined by 
\[
(x,y) \mapsto \mbox{\rm dist}_{D(x)}(y) \qquad (x,y \in X) 
\]
are approximators in the sense of Definition \ref{approximator}.
\end{ass} 

Before proceeding, we make some observations.  Implicit in the assumption (as a consequence of the approximator properties) is that the feasible region $X$ is closed and the objective $f$ is locally Lipschitz.  However, we emphasize that neither need be geodesically convex.  The assumption furthermore ensures that the approximation mapping $D$ has closed graph, and in fact is continuous, since the function $(x,y) \mapsto \mbox{\rm dist}_{D(x)}(y)$ is continuous.  We also see that each approximating set $D(x)$ is a closed second-order approximation of $X$, and each model $h_x$ is a continuous second-order approximation of $f$:  for points $y \in X$ near $x$ we have
\[
\mbox{dist}_{D(x)}(y) = O\big(d^2(x,y)\big) \qquad \mbox{and} \qquad (h_x - f)(y) ~=~ O\big(d^2(x,y)\big).
\]

\begin{exa}[Proximal Point Method]
{\rm
Consider a proper CAT(0) space $X$ and a geodesically convex, locally Lipschitz function $f \colon X \to \R$ that is bounded below.  Let $D(x) = X$ for all points $x \in X$, and for any constant $\mu > 0$, define the associated {\em quadratic penalty modeling function} 
\bmye \label{quadratic-penalty}
h(x,y) ~=~ f(y) + \frac{\mu}{2}d^2(x,y).
\emye
Then Assumption \ref{overall} holds.  The associated envelope is the Moreau-Yosida approximation of $f$, the iteration map is the metric version of the proximal point iteration \cite{bacak}, and MM sequences are just proximal point sequences.
}
\end{exa}

\noindent 
In any proper CAT(0) space $X$, a nontrivial example of an approximating mapping is defined, for any radius $\rho>0$, by $D(x) = B_{\rho}(x)$, as seen in Example \ref{threshhold}.


MM sequences are just fixed point iterations associated with the iteration map $p$.  We first study the fixed point property $x \in p(x)$. Clearly, every minimizer of the objective $f$ is a fixed point of $p$. Moreover, we have the following converse result.

\begin{prop}[Fixed points and zero slope]~ \label{zero}
If Assumption \ref{overall} holds, then at any fixed point of the iteration map, the objective has slope zero.
\end{prop}

\pf
Consider any fixed point $x$.  If the objective $f$ has strictly positive slope, then there exists a constant $\sigma > 0$ and a sequence of points $x_k \ne x$ converging to $x$ in $X$ and all satisfying
\[
\frac{f(x) - f(x_k)}{d(x,x_k)} ~>~ \sigma.
\]
By our assumptions, there exist points $y_k \in D(x)$ satisfying 
$d(y_k,x_k) = o\big(d(x_k,x)\big)$.  The fixed point property ensures $f(x) \le h_x(y_k)$, so for all large $k$ we have
\[
\frac{h_x(x_k) - f(x_k)}{d(x,x_k)} ~=~ \frac{h_x(y_k) - f(x_k)}{d(x,x_k)} ~>~ \sigma,
\]
using Lipschitz property of the model $h_x$.  But this contradicts the implication from our assumptions that $h_x(x_k) - f(x_k) = o\big(d(x_k,x)\big)$.
\finpf

\begin{prop}
On a proper geodesic space, if a geodesically strongly convex function is continuous and bounded below, then it has a unique minimizer.
\end{prop}

\pf
Denote the given function by $g$, and the strong convexity constant by $\mu$.  Fix any point 
$\bar x \in X$.  At any point $x \in X$ satisfying $g(x) \le g(\bar x)$, consider a geodesic with endpoints $x$ and $\bar x$.  The midpoint $x'$ satisfies
\[
\inf g ~\le~ g(x') ~\le~ \frac{1}{2}g(x) + \frac{1}{2}g(\bar x) - \frac{\mu}{8}d^2(x,\bar x)
~\le~ g(\bar x) - \frac{\mu}{8}d^2(x,\bar x),
\]
so $x$ lies in the closed ball defined by $d^2(x,\bar x) \le \frac{8}{\mu}(g(\bar x) - \inf g)$.  By continuity, $g$ has a minimizer $\hat x$ over this compact set, and $\hat x$ must minimize $g$ globally.  The preceding inequality shows uniqueness.
\finpf

\noindent
This result holds more generally if the space is CAT(0) \cite[Lemma 1.3]{mayer}.

\begin{prop}
Suppose that Assumption \ref{overall} holds.  Then the iteration map is single-valued and continuous.  Furthermore, for any associated MM sequence that is bounded, there exists a constant $\beta \ge 0$ such that the envelope $F$ at all successive triples of iterates $x_-, x, x_+$ satisfies the descent conditions
\begin{myeqnarray} 
F(x_-) - F(x) & \ge & \frac{\mu}{2}d^2(x,x_+), \label{envelope_descent} \\
|\nabla F|(x) & \le & \beta d(x,x_+). \label{near-critical}
\end{myeqnarray}
\end{prop}

\pf
The previous result shows that the iteration map $p$ is single-valued.  We next prove the quadratic growth condition
\bmye \label{growth}
h_x(y) - h_x\big(p(x)\big) ~\ge~ \frac{\mu}{2} d^2\big(y,p(x)\big) \qquad \mbox{whenever}~y \in D(x).
\emye 
To see this, for any point $y \ne p(x)$ in the geodesically convex set $D(x)$, consider the distance $\delta = d(p(x),y) > 0$ and a geodesic $\gamma \colon [0,\delta] \to D(x)$  satisfying $\gamma(0) = p(x)$ and $\gamma(\delta) = y$. By 
$\mu$-strong convexity, for $\tau \in [0,\delta]$ we have $\gamma(t) \in D(x)$ and
\[
h_x \big(\gamma(\tau)\big) ~-~ \frac{\mu}{2}\tau^2 ~\le~ \Big(1-\frac{\tau}{\delta}\Big)h_x \big(p(x)\big) ~+~ 
\frac{\tau}{\delta} \Big(h_x(y) - \frac{\mu}{2}\delta^2\Big).
\]
Since $p(x)$ minimizes $h_x$ on $D(x)$, for small $\tau > 0$ we deduce
\[
0 ~\le~ h_x \big(\gamma(\tau)\big) ~-~ h_x \big(p(x)\big) ~\le~
\frac{\tau}{\delta} \Big(h_x(y) - h_x \big(p(x)\big) - \frac{\mu}{2}\delta^2 \Big) ~+~ \frac{\mu}{2}\tau^2.
\]
Dividing by $\tau$ and taking the limit as $\tau \downarrow 0$ proves the growth condition (\ref{growth}).

Consider any sequence $(x_k)$ in $X$ converging to $\bar x$. Since $p(\bar x) \in D(\bar x)$, the continuity of $D$ ensures the existence of points $y_k \in D(x_k)$ satisfying $y_k \to p(\bar x)$. The growth condition (\ref{growth}) implies, for all $k$, the inequality
\[
h_{x_k} (y_k) - h_{x_k} \big(p(x_k)\big) ~\ge~ \frac{\mu}{2} d^2 \big( y_k, p(x_k) \big).
\]
Taking the $\limsup$ of both sides implies $d\big(y_k,p(x_k)\big) \to 0$, and hence
$p(x_k) \to p(\bar x)$.  Thus $p$ is continuous.  

At any point $x \in X$, inequality (\ref{growth}) implies
\begin{eqnarray*}
F(x) - F\big(p(x)\big) 
&=& 
h_x\big(p(x)\big) - h_{p(x)}\Big(p\big(p(x)\big)\Big)  \\
&\ge& h_{p(x)}\big(p(x)\big) - h_{p(x)}\Big(p\big(p(x)\big)\Big) 
~\ge~ \frac{\mu}{2} d^2\Big(p(x),p\big(p(x)\big)\Big).
\end{eqnarray*}
This proves the descent condition (\ref{envelope_descent}).

To prove the final descent condition (\ref{near-critical}), we can assume $|\nabla F|(x) > 0$, since otherwise there is nothing to prove, so there exists a sequence of points $x_k \to x$ in $X$ with
\[
|\nabla F|(x) ~ = ~ \lim_k \frac{F(x) - F(x_k)}{d(x,x_k)} 
~ = ~ \lim_k \frac{h\big(x,p(x)\big) - h\big(x_k,p(x_k)\big)}{d(x,x_k)}.
\]
Consider the sequence
\[
\alpha_k ~=~ d\big(x,p(x_k)\big) + d\big(x_k,p(x_k)\big).
\]
Since $p(x_k) \in D(x_k)$ and the function $(x,y) \mapsto \mbox{dist}_{D(x)}(y)$ is an approximator, we know for some constant $\gamma \ge 0$ the inequality
\[
\mbox{dist}_{D(x)}\big(p(x_k)\big) ~\le~ 
\gamma \alpha_k d(x,x_k)
\]
for all $k$, so there exists a bounded sequence of points $y_k \in D(x)$ satisfying
\[
d\big(y_k,p(x_k)\big) ~\le~ 
\gamma \alpha_k d(x,x_k).
\]
Increasing $\gamma$ if necessary, we can also ensure, for all $k$ the inequalities
\begin{eqnarray*}
h\big(x,p(x_k)\big) - h\big(x_k,p(x_k)\big)
& \le &  \gamma \alpha_k d(x,x_k) \\
h(x,y_k) - h\big(x,p(x_k)\big) & \le & \gamma d\big(y_k,p(x_k)\big),
\end{eqnarray*}
since $h$ is also an approximator.
Using our assumptions and the continuity of $p$, we deduce
\begin{eqnarray*}
|\nabla F|(x)  
& \le & \liminf_k \frac{h\big(x,y_k\big) - h\big(x_k,p(x_k)\big)}{d(x,x_k)} \\
& \le & \liminf_k \frac{h\big(x,p(x_k)\big) + \gamma d\big(y_k,p(x_k)\big) - h\big(x_k,p(x_k)\big)}{d(x,x_k)} \\
& \le & \gamma(\gamma+1) \liminf_k  \alpha_k ~ = ~ 2\gamma(\gamma+1)  d\big(x,p(x)\big).
\end{eqnarray*}
This completes the proof.
\finpf

\noindent
Following the general framework promoted in \cite{taylor}, this result serves to justify the natural ``small step'' stopping criterion for MM algorithms, since if two successive iterates $x,x_+$ are close, then the slope of the envelope at $x$ must be small.

\begin{thm}[MM sequence convergence]~
Suppose that Assumption \ref{overall} holds.  If the corresponding envelope, when restricted to closed balls, always satisfies the KL property for minimization, with a concave desingularizer, then any bounded $MM$ sequence 
$(x_k)$ with initial value $f(x_0)$ near $\inf f$ converges to a minimizer of~$f$.
\end{thm} 

\pf
The proof follows similar techniques in \cite{noll-kl,bolte-pauwels-mm} in imitating that of Theorem~\ref{iterate}.  By our previous results, successive triples of iterates $\ldots x_-,x,x_+ \ldots$ satisfy
\bmye \label{triple}
F(x_-) - F(x) ~ \ge ~ \frac{\mu}{2} d^2(x,x_+) \qquad \mbox{and} \qquad 
d(x,x_+) ~ \ge ~ \beta |\nabla F|(x).  
\emye
We can assume that no iterate minimizes $F$, since otherwise the first inequality would ensure that the sequence is constant thereafter.  Define the gap function $g = F - \inf F$.  The envelope values $F(x_k)$ decrease monotonically, and the KL inequality ensures that the gap $g$ satisfies
\[
1 ~\le~ |\nabla(\phi \circ g)|(x_k)  ~=~ \phi'\big(g(x_k)\big) \cdot |\nabla g|(x_k).
\]
We deduce that each slope $|\nabla g|(x_k) = |\nabla F|(x_k)$ is nonzero, and so the distances 
$d_k = d(x_{k-1},x_k)$ are nonzero for all $k>0$.

Using the concavity of the desingularizer $\phi$, we observe, for all $k>0$, 
\[
\phi\big(g(x_{k-1})\big) - \phi\big(g(x_k)\big) 
\: \ge \: 
\phi'\big(g(x_{k-1})\big) \big(g(x_{k-1}) - g(x_k)\big) 
\: \ge \: 
\frac{F(x_{k-1}) - F(x_k)}{|\nabla F| (x_{k-1})} \: \ge \: 
\frac{\mu d_{k+1}^2}{2\beta d_k},
\]
so
\[
\frac{2\beta}{\mu}\Big( \phi\big(g(x_{k-1})\big) - \phi\big(g(x_k)\big) \Big)
~ \ge ~
\frac{2d_k d_{k+1} - d_k^2}{d_k}
~ \ge ~
2d_{k+1} - d_k.
\]
We deduce that the positive numbers $\lambda_k = \frac{2\beta}{\mu}\phi\big(g(x_{k-1})\big) + d_k$ satisfy the key inequality \mbox{$d(x_k,x_{k+1}) \le \lambda_k - \lambda_{k+1}$}, so the sequence $(x_k)$, being Cauchy, converges to some limit $x^*$.  Theorem \ref{value} shows that $(x_k)$ is a minimizing sequence for $F$ and hence also for $f$, so by continuity we deduce that $x^*$ minimizes $f$, as required.
\finpf

\begin{rem}
{\rm
This result implies the convergence of MM sequences for semi-alge\-braic optimization in broad generality.  More precisely, suppose that the feasible region is a semi-algebraic subset of a Euclidean space, and that the approximating mapping and modeling function are also semi-algebraic.  In that case the envelope $F$ is also semi-algebraic, and so satisfies the condition.
}
\end{rem}

\begin{thm}[Global convergence of MM sequences]~
Suppose that Assumption \ref{overall} holds.  If the corresponding envelope is a KL function with a concave desingularizer, as holds in particular if the feasible region, approximating mapping and modeling function are all semi-algebraic, then any bounded $MM$ sequence converges to a point where the objective has zero slope.
\end{thm}

\pf
We follow the proof of Theorem \ref{global}, but instead using the inequalities (\ref{triple}) as in the previous proof.  If the objective $f$ has bounded level sets, then any MM sequence must be bounded because the corresponding objective values decrease monotonically.  The proof shows that the sequence generated by the iteration $x \leftarrow p(x)$ converges to some limit $x^*$.  Since the iteration map $p$ is continuous, $x^*$ must be a fixed point, from which we deduce $|\nabla f|(x^*) = 0$ by Proposition \ref{zero}.
\finpf

\noindent
Note that MM sequences are always bounded if the objective has bounded level sets.  We end with two illustrations from the literature.

\subsection*{Convex-composite optimization \cite{prx_lin}}
Given Euclidean spaces $X$ and $W$, a ${\mathcal C}^{(2)}$-smooth map $G \colon X \to W$, a convex function $q \colon W \to \R$ and any constant $\mu > 0$, we consider the composite objective $f \colon X \to \R$, a function $H \colon X^2 \to W$ describing the linear approximation to $G$, and the ``prox-linearized'' function 
$h \colon X^2 \to \R$, defined for points $x,y$ in $X$ by
\begin{eqnarray*}
f(x) &=& q\big(G(x)\big),  \\
H(x,y) &=& G(x) + DG(x)[y-x], \\
h(x,y) &=& q\big(H(x,y)\big) ~+~ \frac{\mu}{2} |x-y|^2,
\end{eqnarray*}
where the linear map $DG(x) \colon X \to W$ is the derivative for $G$ at $x$.  With the trivial approximating mapping defined by $D(x) = X$ for all points $x \in X$, we show under reasonable conditions that Assumption \ref{overall} holds.

We suppose that $f$ is bounded below.  We assume too that $q$ is \mbox{$L$-Lipschitz}, and that the second derivative of $G$, the bilinear map $D^2 G(x) \colon X^2 \to W$, has norm bounded by a constant $M$ for all $x$.  (For simplicity of exposition, these are global assumptions;  a more refined development allows the constants to depend on a given bounded domain.)  For any point $x \in X$ and unit vector $w \in W$, consider the function $e \colon X \to \R$ defined by
\[
e(y) ~=~ \ip{w}{H(x,y) - G(y)} \qquad (y \in X).
\]
Notice for all vectors $u \in X$ we have
\[
De(y)[u] = \ip{w}{\big(DG(x) - DG(y)\big)[u]}.
\]
By the Mean Value Theorem, there exists a point $z \in [x,y]$ such that
\[
e(y) ~=~ e(x) + De(z)[y-x] ~=~ \ip{w}{\big(DG(x) - DG(z)\big)[y-x]},
\]
and furthermore a point $u \in [x,z]$ such that 
\[
\big(DG(x) - DG(z)\big)[y-x] ~=~ D^2G(u)[y-x,x-z].
\]
Hence
\[
\big|G(x) + DG(x)[y-x] - G(y)\big| ~\le~  M|y-x|\,|x-z|,
\]
so
\[
\Big|q\big(G(y)\big) ~-~ q\big(G(x) + DG(x)[y-x]\big)\Big| ~\le~  LM|x-y|^2.
\]
Providing $\mu > 2LM$ we deduce the majorization property.

Turning to the remaining assumptions, fix a point $y \in X$.  We first prove that all points $x,z \in X$ satisfy
\bmye \label{H}
|H(x,y) - H(z,y)| ~\le~ M |x-z|\big(|x-y| + |z-y|\big).
\emye
Fix any unit vector $w \in W$, and define a function $c \colon X \to \R$ by 
$c(x) = \ip{w}{H(x,y)}$.  Then, by the Mean Value Theorem, there exists a point $u \in [x,z]$ such that
\[
c(x) - c(z) ~=~ Dc(u)[x-z] ~=~ \ip{w}{D^2G(u)[y-u,x-z]}.
\]
Since $w$ was arbitrary, our claim now follows, since $|y-u| \le |x-y| + |z-y|$.

Using inequality (\ref{H}), we deduce
\begin{eqnarray*}
|h(x,y) - h(z,y)| 
&\le&
LM |x-z|\big(|x-y| + |z-y|\big) ~+~ 
\frac{\mu}{2}\big| |x-y|^2 - |z-y|^2\big| \\
&=& 
LM |x-z|\big(|x-y| + |z-y|\big) ~+~ 
\mu |\ip{x-z}{x+z-2y}| \\
& \le & 
\Big(LM + \frac{\mu}{2}\Big) |x-z|\big(|x-y| + |z-y|\big).
\end{eqnarray*}
Assumption \ref{overall} follows.

Suppose in addition that the functions $q$ and $G$ are semi-algebraic.  Then for all sufficiently large constants $\mu > 0$, any bounded sequence generated by the prox-linear iteration 
\[
x ~\leftarrow~ \mbox{arg}\!\min_y \Big\{ q\big(G(x) + DG(x)[y-x]\big) ~+~ \frac{\mu}{2}|x-y|^2 \Big\}
\]
converges to a point where the composite function $q \circ G$ has slope zero.

\subsection*{Moving balls method \cite{moving_balls,bolte-pauwels-mm}}
For differentiable functions  $f, f_1,f_2, \ldots, f_m \colon \R^n \to \R$ with Lipschitz gradients, we consider the problem of minimizing the objective $f$ over the feasible region
\[
X ~=~ \{x \in \R^n : f_i (x) \le 0~ \mbox{for each}~i\}.
\] 
We assume that this classical nonlinear program is bounded below, and that the Mangasarian-Fromovitz constraint qualification (MFCQ) holds:  in other words, for all points \mbox{$x \in X$}, there exists a vector $d \in \R^n$ satisfying
\[
\langle \nabla f_i(x), d\rangle < 0 \qquad \mbox{for all}~ i \in I(x) =  \{i \colon f_i(x) = 0\}.
\]

With the aim of constructing an optimization method that maintains feasibility, our majorization-minimization tools consist of an approximating mapping $D$ and a modeling function $h$ defined by
\begin{eqnarray*}
h(x,y) &=& f(x) + \langle \nabla f(x), y-x\rangle + \frac{L}{2}\|y-x\|^2, \\
h_i(x,y) &=& f_i(x) + \langle \nabla f_i(x), y-x\rangle + \frac{L_i}{2}\|y-x\|^2 \quad \mbox{for each}~i, \\
D(x) & = & \bigcap_i D_i(x) \quad \mbox{with}~~D_i(x) = \{ y \in \R^n : h_i(x,y) \le 0 \} ~~ \mbox{for each}~i,
\end{eqnarray*}
where $L, L_1, \dots, L_m$ are the Lipschitz moduli of $\nabla f, \nabla f_1, \dots, \nabla f_m$ respectively.
For each point $x \in X$, the set $D_i(x)$ is a closed ball centered at the point 
$x - \frac{1}{L_i}\nabla f_i(x)$ and with radius
\[
\rho_i(x) ~=~ \sqrt{\frac{1}{L_i^2} \|\nabla f_i(x)\|^2 - \frac{2}{L_i} f_i (x)}.
\]
Hence the name of the method.  Our aim is to show that Assumption \ref{overall} holds, ensuring convergence of the method under reasonable conditions.

It suffices to verify that both the functions $h$ and $(x,y) \mapsto \mbox{\rm dist}_{D(x)}(y)$ (for points \mbox{$x,y \in \R^n$}) are approximators, since the other conditions easily follow from the standard quadratic upper bound for $L$-smooth functions $g \colon \R^n \to \R$:
\[
g(y) ~\le~ g(x) + \langle \nabla g(x), y-x\rangle + \frac{L}{2}\|y-x\|^2.
\]
The fact that $h$ is an approximator follows from Examples \ref{squared} and \ref{linear}, along with Proposition \ref{vector}.  In fact, the following inequality is easy to verify directly that all points $x,y,z \in \R^n$ satisfy the inequality
\bmye \label{analog}
|h (x, y) - h (z, y)| ~\le~ \frac{3L}{2} |x-z|\big(|x-y| + |z-y|\big)
\emye
along with the corresponding inequality for the functions $h_i$.

To show that the function $(x,y) \mapsto \mbox{\rm dist}_{D(x)}(y)$ is also an approximator, we need the following elementary fact \cite{Nedic2010}: for closed convex sets $S_1,S_2, \ldots, S_m \subset \R^n$, each containing the ball $B_{\delta} (0)$, all points $x \in \R^n$ satisfy
\bmye\label{cvx-intersect}
\mbox{\rm dist}_{\cap_i S_i}(x) ~\le~ \frac{\|x\|}{\delta} \max_i \big\{\mbox{\rm dist}_{S_i}(x)\big\}.
\emye

Fix a point $\bar{x} \in X$. By the MFCQ condition, there exists a unit vector $d \in \R^n$ satisfying
\[
\langle \nabla f_i(\bar{x}), d\rangle < 0 \quad \mbox{for all}~ i \in I(\bar{x}). 
\]
For any constant $\epsilon \in (0, \min_i 2L_i)$ sufficiently small, all points $x \in B_{\epsilon}(\bar{x})$ satisfy the inequalities
\bmye\label{bound-fx-nablafx}
f_i(x) \le -\epsilon \quad \mbox{for all}~ i \notin I(\bar{x}) \qquad \text{and} \qquad 
\nabla f_i(x)^{\top} d \le -\epsilon \quad \mbox{for all}~ i \in I(\bar{x}).
\emye
and, for sufficiently large $M \ge \max_i L_i$, the inequalities $\|\nabla f_i(x)\| \le M$ for each $i=1,2,\ldots,m$.

Define constants $t = \min\{1, \frac{\epsilon}{4M}\}$ and $\delta = \min\{1, \frac{\epsilon}{2M}\} t$.  For any point $x \in B_{\epsilon} (\bar{x}) \cap X$ and any vector $e \in B_{\delta} (0)$, consider the quantity
\[
h_i (x, x+td+e) ~=~ f_i (x) + \langle \nabla f_i(x), td+e \rangle + \frac{L_i}{2} \|td+e\|^2.
\]
For all $i \notin I(\bar{x})$, the left-hand side is no larger than $-\epsilon + 4Mt$, and for all 
$i \in I(\bar{x})$ it is no larger than $-\frac{\epsilon}{2}t + 2Mt^2$, so in either case it is nonpositive.  We deduce 
\[
B_{\delta}(x+td) ~\subset~ \bigcap_i D_i(x) ~=~ D(x). 
\]
Property (\ref{cvx-intersect}) therefore implies that all points $y \in \R^n$ satisfy the inequality
\bmye\label{dist-intersect}
\mbox{\rm dist}_{D(x)}(y) ~\le~ \frac{\|y-x-td\|}{\delta} \max_i \big\{ \mbox{\rm dist}_{D_i(x)}(y) \big\}.
\emye
Consider points $x \in B_{\epsilon} (\bar{x}) \cap X$, $z \in X$ and $y \in D(z)$.  If 
$y \notin D_i (x)$, then
\begin{eqnarray*}
\mbox{\rm dist}_{D_i(x)}(y)
&=& \big\| y - x + \frac{1}{L_i} \nabla f_i(x) \big\| ~-~ \rho_i(x) \\
&=& \frac{2}{L_i} \cdot \frac{h_i (x, y)}{\big\| y - x + \frac{1}{L_i} \nabla f_i(x) \big\| ~+~ \rho_i(x)} \\ \\
&\le& \frac{2}{L_i} \cdot \frac{h_i (x, y) - h_i (z, y)}{\big\| y - x + \frac{1}{L_i} \nabla f_i(x) \big\| ~+~ \rho_i(x)} \\ \\
&\le& \frac{h_i (x, y) - h_i (z, y)}{\sqrt{\|\nabla f_i(x)\|^2 - 2 L_i f_i (x)}} 
~\le~ \frac{1}{\epsilon} \big(h_i (x, y) - h_i (z, y)\big) \\
&\le& \frac{3L_i}{2\epsilon} |x-z|\big(|x-y| + |z-y|\big),
\end{eqnarray*}
the penultimate inequality following from the inequalities (\ref{bound-fx-nablafx}).  Property (\ref{dist-intersect}) now implies
\[
\mbox{\rm dist}_{D(x)}(y) ~\le~ \frac{3 M \|y-x-td\|}{2 \epsilon \delta} |x-z|\big(|x-y| + |z-y|\big).
\]
Since the above inequality holds for all $y \in D(z)$, it along with (\ref{dist-approximator}) implies the local Lipschitz continuity of the function $(x,y) \mapsto \mbox{\rm dist}_{D(x)}(y)$.
Hence, by Lemma \ref{lem:approximator}, the function $(x,y) \mapsto \mbox{\rm dist}_{D(x)}(y)$ is an approximator, completing our verification of Assumption \ref{overall}.


\def\cprime{$'$} \def\cprime{$'$}

\end{document}